\documentclass[11pt]{amsart}
\usepackage{microtype}
\usepackage{mathtools}
\usepackage{graphicx}
\usepackage{comment,overpic,epsfig}
\usepackage[usenames,dvipsnames]{color}
\usepackage[hyphens]{url}
\usepackage[pagebackref,linktocpage=true,colorlinks=true,linkcolor=RoyalBlue,citecolor=BrickRed,urlcolor=RoyalBlue]{hyperref}
\usepackage[msc-links,abbrev]{amsrefs}
\usepackage[safe]{tipa}

\newtheorem{thm}{Theorem}[section]
\newtheorem{lem}[thm]{Lemma}
\newtheorem{cor}[thm]{Corollary}
\newtheorem{prop}[thm]{Proposition}

\theoremstyle{definition}

\newtheorem{note}[thm]{Note}

\theoremstyle{remark}

\newcommand{\R}{\mathbf{R}}

\newcommand{\ol}{\overline}
\newcommand{\lan}{\langle}
\newcommand{\ran}{\rangle}

\newcommand{\C}{\mathcal{C}}

\newcommand{\arc}[1]{{%
  \setbox9=\hbox{#1}%
  \ooalign{\resizebox{\wd9}{\height}{\texttoptiebar{\phantom{A}}}\cr#1}}}

\renewcommand{\S}{\mathbf{S}}
\renewcommand{\tilde}{\widetilde}
\renewcommand{\hat}{\widehat}
\renewcommand{\epsilon}{\varepsilon}

\DeclareMathOperator{\dist}{dist}
\DeclareMathOperator{\inte}{int}
\DeclareMathOperator{\conv}{conv}

\DeclareMathOperator{\ave}{ave}
\DeclareMathOperator{\cm}{cm}
\DeclareMathOperator{\Emb}{Emb}
\DeclareMathOperator{\Imm}{Imm}
\DeclareMathOperator{\length}{length}
\DeclareMathOperator{\supp}{supp}

\DeclareMathOperator{\SL}{SL}

\begin{document}


\title{Deformations of Curves with Constant Curvature}

\author{Mohammad Ghomi}
\address{School of Mathematics, Georgia Institute of Technology,
Atlanta, Georgia 30332}
\email{ghomi@math.gatech.edu}
\urladdr{www.math.gatech.edu/$\sim$ghomi}

\author{Matteo Raffaelli}
\address{School of Mathematics, Georgia Institute of Technology, Atlanta, Georgia 30332}
\email{raffaelli@math.gatech.edu}
\urladdr{matteoraffaelli.com}

\date{Last revised on \today}
\thanks{The first named author was supported by NSF grant DMS-2202337.}

\begin{abstract}
We  prove that curves of constant curvature satisfy the parametric $\mathcal{C}^1$-dense relative $h$-principle in the space of  immersed curves with nonvanishing curvature in Euclidean space $\R^{n\geq 3}$.  It follows that two knots of constant curvature in $\R^3$ are isotopic, resp. homotopic,  through curves  of constant curvature if and only if they are  isotopic, resp.\ homotopic, and their self-linking numbers, resp. self-linking numbers mod $2$, are equal. The proofs are based on convex integration techniques, utilizing parametric versions of classical theorems by Carath\'{e}odory and Steinitz in convex geometry.
\end{abstract}

\maketitle

\section{Introduction}
The homotopy and isotopy classes of closed curves with nonvanishing curvature in Euclidean space $\R^n$ have been classified in terms of their self-linking number \cite{feldman:curve,rosenberg:constant,gluck&pan}. In this work we obtain parallel results for curves of constant curvature by an approximation argument. 

To state our main result, let $\Gamma$ be an interval $[a,b]\subset\R$  or the topological circle $\R/((b-a)\mathbf{Z})$, and $\mathcal{C}^k(\Gamma,\R^n)$ be the space of $\C^k$ maps $f\colon \Gamma\to\R^n$ with its standard topology induced by the $\C^k$-norm 
$|\cdot|_k$.  We say that $f_0$, $f_1\in\C^k(\Gamma,\R^n)$ are \emph{$\C^{\ell\leq k}$-homotopic} if there is a family  $f_t\in\C^\ell(\Gamma,\R^n)$, for $t\in[0,1]$, such that $t\mapsto f_t$ is continuous. Then $f_t$ is called a $\C^\ell$-\emph{homotopy}. The space of (immersed) \emph{curves} $\Imm^{k}(\Gamma,\R^n)\subset \mathcal{C}^k(\Gamma,\R^n)$ consists of locally injective maps with nonvanishing derivative when $k\geq 1$.  
 The \emph{curvature} of $f\in\Imm^{2}(\Gamma,\R^n)$ is given by
$
\kappa\coloneqq |T'|/|f'|,
$
where
$
T\coloneqq f'/|f'|
$
is the \emph{tantrix} of $f$.

\begin{thm}\label{thm:main1}
Let $f_t\in \Imm^{k\geq 2}(\Gamma,\R^{n\geq 3})$ be a $\C^2$-homotopy of curves with nonvanishing curvature. Suppose that
$f_0$ and $f_1$ have constant curvature. Then, for any $\epsilon>0$, there exists a $\C^2$-homotopy $\tilde f_t\in\Imm^k(\Gamma,\R^n)$  of curves with constant curvature such that  $\tilde f_0=f_0$,
$\tilde f_1=f_1$, and $|\tilde f_t- f_t|_1\leq\epsilon$.
\end{thm}

In the terminology of Gromov \cite{gromov:PDR} or Eliashberg--Mishachev \cite{eliashberg&mishachev}, the result above establishes  a parametric $\C^1$-dense relative $h$-principle for curves of constant curvature in the space of curves with nonvanishing curvature. The nonparametric version of Theorem 
\ref{thm:main1}, which states that curves of constant curvature are dense in $\Imm^{k\geq2}(\Gamma,\R^n)$  with respect to the $\C^1$-norm, 
had been obtained earlier in \cite{ghomi:knots}, see also \cite{ghomi-raffaelli2025}.  A result analogous to Theorem \ref{thm:main1} is Bleecker's parametric version \cite{bleecker1997} of the 
Nash--Kuiper theorem \cite{nash1954,kuiper1955}: a continuous family of 
short embeddings of a compact Riemannian $n$-manifold into $\R^{n+1}$ may be approximated in the $\C^0$-norm by a family of $\C^1$ isometric embeddings; however, this result does not fix the endpoints if they are already isometric. See \cite{ghomi&kossowski} for another example of the parametric $h$-principle in Riemannian geometry.

We say that $f\in\C^k(\Gamma,\R^n)$  is \emph{closed} when $\Gamma$ is a circle. 
If closed curves $f_0$, $f_1\in\Imm^{k\geq 2}(\Gamma,\R^{3})$ with nonvanishing curvature are $\C^2$-homotopic through curves with nonvanishing curvature, then their tantrices $T_0$, $T_1$ are $\C^1$-homotopic through immersed curves. Feldman \cite{feldman:curve} proved the converse of this phenomenon, via Smale's extension \cite{smale:homotopy} of the Whitney--Graustein theorem \cite{whitney1937} to spherical curves. Furthermore, Rosenberg \cite{rosenberg:constant} showed that $T_0$ and $T_1$ are $\C^1$-homotopic through immersed curves if and only if $\SL_2(f_0)=\SL_2(f_1)$, where $\SL_2$ denotes the self-linking number mod $2$. Thus Theorem \ref{thm:main1} yields:

\begin{cor}
A pair of closed curves $f_0$, $f_1\in\Imm^{k\geq 2}(\Gamma,\R^{3})$ of constant curvature are $\C^2$-homotopic through $\C^k$ curves of constant curvature  if and only if $\SL_2(f_0)=\SL_2(f_1)$.
\end{cor}

Let $\Emb^k(\Gamma,\R^n)\subset\Imm^k(\Gamma,\R^n)$ denote the space of injective or embedded curves. 
A pair of curves $f_0$, $f_1\in \Emb^k(\Gamma,\R^n)$ are $\C^k$-\emph{isotopic} if they can be joined by a $\C^k$-homotopy $f_t\in \Emb^k(\Gamma,\R^n)$. A \emph{knot} is an embedded closed curve $f\in\Emb^0(\Gamma,\R^3)$.  Gluck and Pan \cite{gluck&pan} showed that two knots of nonvanishing curvature $f_0$, $f_1\in \Emb^{k\geq 2}(\Gamma,\R^3)$ are $\C^{2}$-isotopic through knots of nonvanishing curvature if and only if they are $\C^0$-isotopic and $\SL(f_0)=\SL(f_1)$, i.e., they have the same self-linking number.  If $f_t$ in Theorem \ref{thm:main1} is an isotopy, then (choosing $\epsilon$ sufficiently small) we may assume that $\tilde f_t$ is an isotopy as well. Thus we obtain:

\begin{cor}
A pair of knots $f_0$, $f_1\in\Emb^{k\geq 2}(\Gamma,\R^{3})$ of constant curvature are $\C^2$-isotopic through $\C^k$ knots  of constant curvature if and only if they are $\C^0$-isotopic and $\SL(f_0)=\SL(f_1)$.
\end{cor}

So the $\C^2$ homotopy and isotopy classes of closed curves of constant curvature in $\R^3$ mirror those of closed curves with nonvanishing curvature. In other words, a loop of constant curvature is just as flexible as one without inflection points. The self-linking number of a knot $f\in\Emb^2(\Gamma,\R^3)$ with nonvanishing curvature is the linking number of $f$ with a perturbation of $f$ along its principal normal $N:=T'/|T'|$. The self-linking number mod $2$ of a closed curve $f\in\Imm^2(\Gamma,\R^3)$ with nonvanishing curvature is the self-linking number mod $2$ of a perturbation of $f$ to an embedded curve. See \cite{gluck&pan, rosenberg:constant,pohl1968} for further background on self-linking number.

The proof of Theorem~\ref{thm:main1} utilizes the fact that, when $|f_t'|=1$, the curvature of $f_t$ is the speed of its tantrix $T_t$. Thus we approximate $T_t$ with a homotopy $\tilde T_t$ of curves with constant speed, and set $\tilde f_t(s):=\int_a^s\tilde T_t(u)du$. For $\tilde f_t$ to be closed when $f_t$ is closed, we need to have $\int_\Gamma\tilde T_t=\int_\Gamma T_t$. To this end, we construct $\tilde T_t$ by creating sinusoidal bumps near suitably chosen points along $T_t$, via parametric versions of classical theorems of Carath\'{e}odory and Steinitz, which may be of independent interest. The bumps have to be long enough to control $\int_\Gamma\tilde T_t$, while vanishing sufficiently fast as $t\to 0$, $1$ to ensure that $\tilde T_t$ is a $\C^1$-homotopy. Implementing this plan is particularly delicate when the convex hull of $T_0$ 
or $T_1$ has no interior points, i.e., when these tantrices are \emph{flat}.
 
Constructing closed curves of constant curvature, by integrating spherical curves, goes back to Fenchel \cite{fenchel1930}, and  belongs to the convex integration theory \cite{gromov:PDR,eliashberg&mishachev,spring1998,geiges2003}, which includes works by Whitney \cite{whitney1937} and Nash--Kuiper  \cite{nash1954,kuiper1955}.
 A nonparametric version of the above argument had been developed in \cite{ghomi:knots}; however, the perturbation method there is not suitable for obtaining a $\C^1$-homotopy of tantrices. Our present construction is fundamentally different, and more demanding technically.

This paper is organized as follows. In Section \ref{sec:prelim} we show, via a jet
transversality argument, that the tantrices of the initial homotopy in Theorem
\ref{thm:main1} may be assumed to be nonflat for $t\in(0,1)$. In Section
\ref{sec:geod} we construct perturbations of circular arcs on a sphere which increase their
length while preserving their center of mass. Section \ref{sec:lemmas} develops
parametric versions of the convex hull theorems of Steinitz and Carath\'eodory. These results
are combined in Section \ref{sec:spherical} to prescribe the length and center of
mass of homotopies of nonflat spherical curves. Finally, in Section
\ref{sec:proof} we apply this result to tantrices and prove Theorem
\ref{thm:main1}.

Nontrivial $\mathcal{C}^2$ knots of constant curvature were first constructed by Koch and Engelhardt \cite{koch:constant} by gluing helical segments. Using the same method, McAtee \cite{mcatee2007}  obtained $\mathcal{C}^2$ knots of constant curvature in every isotopy class.  Existence of $\C^\infty$ knots of constant curvature, in every isotopy class, was first established in \cite{ghomi:knots} via convex integration. See also \cite{ghomi-raffaelli2024torsion} for a quicker proof with the aid of degree theory. For knots with prescribed curvature, an $h$-principle was first obtained by Wasem \cite{wasem2016} in the $\C^2$ category, which was recently extended to the $\C^\infty$ case \cite{ghomi-raffaelli2025}.

 \section{Approximation by Curves with Nonflat Tantrices}\label{sec:prelim}
We always assume that $t\in[0,1]$ unless specified otherwise, and set $I:=[a,b]$.  Here we show  that $\C^2$-homotopies $f_t\in\Imm^{k}(I,\R^{n})$ with nonvanishing curvature may be perturbed to assume  a generic property. A subset of $\R^n$ is \emph{nonflat} if its convex hull has interior points. By abuse of notation, for any mapping $f\in\C^k(I,\R^n)$, we may write $f$ to refer  to $f(I)$, when the meaning is clear from  context, e.g., we say $f$ is nonflat if $f(I)$ is nonflat. First we record the following fact, which is a standard application of Thom's jet transversality theorem \cite{hirsch:book,eliashberg&mishachev}.

\begin{lem}\label{lem:thom}
Let $f\in\C^{k\geq 1}( I ,\R^{n\geq 3})$. Then for any $\epsilon>0$ there exists $\ol f\in \Imm^{\infty}( I ,\R^n)$ such that $|\ol f-f|_k\leq\epsilon$ and the tantrix $\ol T$ of $\ol f$ is nonflat on every subinterval of $I$.
\end{lem}
\begin{proof}
Let $J^{n+1}(I,\R^n)=I\times\R^{n(n+2)}$ be the space of $(n+1)$-jets of maps 
$f\in\mathcal{C}^{n+1}(I,\R^n)$, and $J^{n+1}_*\subset J^{n+1}(I,\R^n)$ be the 
open subset where $f'\neq 0$. For any curve $f$ with $f'\neq 0$ and tantrix 
$T:=f'/|f'|$, a straightforward induction, using 
$|f'|'=\lan f',f''\ran/|f'|$, shows that
$
T^{(j)}=Q_j/(|f'|^{2j+1}),
$
where the entries of $Q_j$ are polynomials in $f',\dots,f^{(j+1)}$. Thus 
$T',\dots,T^{(n)}$ may be regarded as functions on $J^{n+1}_*$, and 
$$
A:=\big\{\det\big[\,T',\dots,T^{(n)}\big]=0\big\}\subset J^{n+1}_*
$$
is the zero set of the polynomial $\det[\,Q_1,\dots,Q_n]$ in the jet 
coordinates. So $A$ is a real algebraic subset of each fiber of $J^{n+1}_*$, 
and thus admits a stratification.

We claim that $A$ has codimension $\geq 1$ in $J^{n+1}_*$. To this end, it suffices to check that $\det[\,Q_1,\dots,Q_n]$ does not vanish 
identically on any fiber of $J^{n+1}_*$, for then $A$ meets each fiber in a 
proper algebraic subset, which has codimension $\geq 1$ in the fiber. To this 
end, it is enough to find, in each fiber, the jet of a unit speed curve $f$ 
with tantrix $T=f'$ such that $T'(s_0),\dots,T^{(n)}(s_0)$ are linearly 
independent; note that 
$$
\det[\,T',\dots,T^{(n)}]=\det[\,Q_1,\dots,Q_n]
$$ 
for 
such jets, since $|f'|\equiv 1$. The desired jet may be generated as follows. 
Curves in $\S^{n-1}$ passing through a point $p$ at time $s_0$ may be written 
as 
$
T=\sqrt{1-|w|^2}\,p+w,
$
where $w$ is a curve in $p^\perp$ with $w(s_0)=0$. Then 
$T^{(j)}(s_0)=w^{(j)}(s_0)+\alpha_j p$, where $\alpha_j$ depends only on 
derivatives of $w$ of order less than $j$, and $\alpha_2=-|w'(s_0)|^2$. So 
setting $w'(s_0):=e_1$, $w''(s_0):=0$, and $w^{(j)}(s_0):=e_{j-1}$ for 
$3\leq j\leq n$, where $e_i$ is an orthonormal basis of $p^\perp$, yields 
that $T'(s_0),\dots,T^{(n)}(s_0)$ are given by $e_1$, $-p$, 
$e_2+\alpha_3 p,\dots, e_{n-1}+\alpha_n p$, which are linearly independent. 
This establishes the claim.

First approximate $f$ in the $\C^k$ norm by a smooth map. By the jet
transversality theorem, this approximation may be chosen arbitrarily
$\C^k$-close to $f$ so that its $(n+1)$-jet is transverse to every
stratum of $A$. Since immersions form an open dense subset of
$\C^\infty(I,\R^n)$, we may also assume that the resulting map
$\ol f$ is immersed. Since $\dim(I)=1\leq\textup{codim}(A)$, transversality implies that
$j^{n+1}\ol f$ meets $A$ only at a discrete set of points. Finally, suppose that the restriction of $\ol T$ to a subinterval of $I$ 
lies in a hyperplane $\lan x,\nu\ran=c$, $\nu\neq 0$. Then 
$\lan\ol T^{(j)},\nu\ran\equiv 0$ on the corresponding segment $J\subset I$, 
for $1\leq j\leq n$. So $\ol T',\dots,\ol T^{(n)}$ lie in the hyperplane 
$\nu^\perp$, and consequently $\det[\,\ol T',\dots,\ol T^{(n)}]\equiv 0$ on 
$J$, which is impossible.
\end{proof}

Using the above observation we now establish the main result of this section:

\begin{prop}\label{prop:nonflat}
Let $f_t\in\Imm^{k\geq 2}(I,\R^{n\geq 3})$ be a
$\C^2$-homotopy with nonvanishing curvature. Then, for
any $\epsilon>0$, there exists a $\C^2$-homotopy
$\ol f_t\in\Imm^k(I,\R^n)$ with nonvanishing curvature
such that $\ol f_0=f_0$, $\ol f_1=f_1$,
$|f_t-\ol f_t|_2\leq\epsilon$, the tantrix $\ol T_t$
of $\ol f_t$ is nonflat for $t\in(0,1)$,
$\ol f_t=f_t$ on a neighborhood of $\partial I$, and
$\length(\ol f_t)=\length(f_t)$.
\end{prop}
\begin{proof}
For $i\in\mathbf{Z}$ choose $t_i\in(0,1)$ with
$t_i<t_{i+1}$ and $t_i\to1$ or $0$ as $i\to\infty$
or $-\infty$, respectively. By Lemma \ref{lem:thom},
there are curves $g_{t_i}\in\Imm^\infty(I,\R^n)$ such
that the tantrix of $g_{t_i}$ is, on every subinterval,
nonflat, and
$
|g_{t_i}-f_{t_i}|_2\to0
$
as $i\to\pm\infty$. For $\lambda\in[0,1]$, set
\begin{equation}\label{eq:olf}
\ol f_{(1-\lambda)t_i+\lambda t_{i+1}}
:=
(1-\lambda)g_{t_i}+\lambda g_{t_{i+1}}.
\end{equation}
This defines a homotopy $\ol f_t\in\C^\infty(I,\R^n)$
for $t\in(0,1)$ with
$|\ol f_t-f_t|_2\to0$ as $t\to0$, $1$. Thus, setting
$\ol f_0:=f_0$ and $\ol f_1:=f_1$ yields a
$\C^2$-homotopy of $\C^k$ curves for $t\in[0,1]$.

Since $t\mapsto f_t$ is continuous with respect to the
$\C^2$-norm, by Lemma \ref{lem:thom} we may choose
$t_i$ and $g_{t_i}$ so that
\begin{equation}\label{eq:ft-ft-prime}
|f_t-f_{t'}|_2\leq\delta
\quad\text{for}\quad
t,t'\in[t_i,t_{i+1}],
\qquad
|g_{t_i}-f_{t_i}|_2\leq\delta,
\end{equation}
for any given $\delta>0$. For
$t\in[t_i,t_{i+1}]$, we may rewrite \eqref{eq:olf} as
$$
\ol f_t
=
(1-\lambda_t)g_{t_i}+\lambda_tg_{t_{i+1}},
$$
where $\lambda_t\in[0,1]$ is determined by
$t=(1-\lambda_t)t_i+\lambda_tt_{i+1}$. Then
\begin{gather}\label{eq:ol-ft-ft}
|\ol f_t-f_t|_2
\leq
(1-\lambda_t)|g_{t_i}-f_t|_2
+
\lambda_t|g_{t_{i+1}}-f_t|_2
\\
\leq
|g_{t_i}-f_{t_i}|_2
+
|f_{t_i}-f_t|_2
+
|g_{t_{i+1}}-f_{t_{i+1}}|_2
+
|f_{t_{i+1}}-f_t|_2
\leq4\delta.
\notag
\end{gather}
Thus, choosing $\delta$ sufficiently small, we obtain
$|\ol f_t-f_t|_2\leq\epsilon$ and ensure that
$\ol f_t$ is immersed and has nonvanishing curvature.

To make the tantrix $\ol T_t$ of $\ol f_t$ nonflat,
let $\phi\colon I\to[0,1]$ be a $\C^\infty$ step
function with $\phi=0$ near $a$ and $\phi=1$ near $b$.
Set $\ol t_i:=(t_i+t_{i+1})/2$, and
$$
g_{\ol t_i}
:=
\phi g_{t_i}+(1-\phi)g_{t_{i+1}}.
$$
Then the tantrix of $g_{\ol t_i}$ is nonflat, since
$g_{\ol t_i}$ coincides with $g_{t_i}$ on a subinterval,
and hence so do their tantrices. We have
$$
|g_{\ol t_i}-f_{\ol t_i}|_2
\leq
C|g_{t_i}-g_{t_{i+1}}|_2
+
|g_{t_{i+1}}-f_{\ol t_i}|_2,
$$
where $C$ depends only on $|\phi|_2$. Furthermore,
$$
|g_{t_i}-g_{t_{i+1}}|_2
\leq
|g_{t_i}-f_{t_i}|_2
+
|f_{t_i}-f_{t_{i+1}}|_2
+
|f_{t_{i+1}}-g_{t_{i+1}}|_2
\leq3\delta,
$$
$$
|g_{t_{i+1}}-f_{\ol t_i}|_2
\leq
|g_{t_{i+1}}-f_{t_{i+1}}|_2
+
|f_{t_{i+1}}-f_{\ol t_i}|_2
\leq2\delta.
$$
Thus
$$
|g_{\ol t_i}-f_{\ol t_i}|_2
\leq
(3C+2)\delta.
$$
Set $\delta\leq\epsilon/(4(3C+2))$, replace
$\{t_i\}$ by a reindexing of $\{t_i,\ol t_i\}$, and
let $\ol f_t$ be the corresponding homotopy given by
\eqref{eq:olf}. Then
$|g_{t_i}-f_{t_i}|_2\leq\epsilon/4$, which yields
$|\ol f_t-f_t|_2\leq\epsilon$ by
\eqref{eq:ol-ft-ft}. Moreover, $\ol T_t$ is
nonflat, since $\ol f_t$ coincides with $g_{t_i}$ or
$g_{t_{i+1}}$ on a subinterval when
$t\in[t_i,t_{i+1}]$. Using a partition of unity on $I$, we may glue
$\ol f_t$ to $f_t$ on sufficiently small neighborhoods
of $\partial I$, after adjusting $\delta$, leaving intact the segments on which
$\ol f_t$ coincides with $g_{t_i}$ or
$g_{t_{i+1}}$.

Finally, we arrange that
$\length(\ol f_t)=\length(f_t)$. Let $J$ be a segment
in the region where $0<\phi<1$; thus $I\setminus J$
contains, for each $t$, a segment on which $\ol f_t$
coincides with one of the curves $g_{t_i}$, and hence
$\ol T_t$ is nonflat on $I\setminus J$. Let
$\psi\geq0$ be a $\C^\infty$ function supported in
$J$ with $\psi>0$ on a subsegment of $J$ about a point
$s_*$, and set 
$$
w_t:=\psi u_t,
\qquad\quad
u_t
:=
\frac{\ol T_t'(s_*)}{|\ol T_t'(s_*)|}.
$$
This is well-defined since $\ol f_t$ has nonvanishing
curvature. By the first variation formula for length,
$$
\left.
\frac{\partial}{\partial\lambda}
\right|_{\lambda=0}
\length\big(\,\ol f_t+\lambda w_t\big)
=
-\int_J\psi\big\lan\,\ol T_t',u_t\big\ran\,ds.
$$
This quantity is negative, and its absolute value is bounded away from
zero uniformly in $t$, provided that $\supp(\psi)$ is sufficiently small,
since $\lan\ol T_t',u_t\ran$ is then uniformly close to
$|\ol T_t'(s_*)|$ on the support of $\psi$, which is bounded below
by a positive constant by the compactness of the family.

Since $\length(\ol f_t)-\length(f_t)$ is continuous,
vanishes at $t=0$, $1$, and can be made uniformly
arbitrarily small by choosing the preceding approximation
sufficiently close, the implicit function theorem gives
a continuous choice of $\lambda_t$, with
$\lambda_0=\lambda_1=0$ and $|\lambda_t|$ arbitrarily
small, such that
$\ol f_t+\lambda_tw_t$ has the same length as $f_t$.
As $w_t$ is supported in $J$ and $|\lambda_t|$ is small, this
adjustment preserves
the properties of $\ol f_t$ established above, and the
preceding approximation may be chosen so that the final
homotopy remains within $\epsilon$ of $f_t$.
\end{proof}

\section{Perturbation of Circular Arcs} \label{sec:geod}
A curve $f\in \Emb^{\infty}(I,\S^{n-1})$ is a \emph{circular arc} if it has constant speed and lies in a plane in $\R^n$. Thus circular arcs generalize geodesics, which lie in planes passing through the origin.
Here we show that a circular arc can be perturbed to a curve of any prescribed greater length which coincides with the arc near its endpoints, and has the same center of mass.
The \emph{center of mass} and \emph{average} of a curve $f\in\Imm^1(I,\R^n)$ are defined as
$$
\cm(f):=\frac{1}{\length(f)}{\int_If|f'|},\quad\quad\text{and}\quad\quad\ave(f):=\frac{1}{|I|}\int_I f.
$$
Note that $\cm(f)$ is invariant under reparametrizations of $f$, and when $|f'|$ is constant,  $\cm(f)=\ave(f)$. Let $\conv(f)$ denote the convex hull of $f$. 
If $f$ is nonflat, then  $\cm(f)$ and $\ave(f)$ lie in the interior of $\conv(f)$ \cite[Lem.\ 2.3]{ghomi:knots} which is denoted by $\inte\conv(f)$. A subset of an affine subspace $V\subset\R^n$ is \emph{nonflat relative to $V$} if its convex hull has interior points in $V$.

\begin{lem}\label{lem:bumps}
For $t\in(0,1)$, let $f_t\in\Emb^{\infty}(I,\S^{n-1})$, $n\geq 3$, be a $\C^{1}$-isotopy of circular arcs of length $\ell_t$, subtending angles at most $\pi/2$ in their circles, and suppose that these circles have a common center. Then, for any continuous family of constants $\tilde\ell_t>\ell_t$ and $\epsilon>0$,  there exists a $\C^1$-isotopy $\tilde f_t\in\Emb^{\infty}(I,\S^{n-1})$ such that $|\tilde f_t-f_t|_0\leq\epsilon$, $\tilde f_t=f_t$ on a neighborhood of $\partial I$ of radius $|I|/32$, $\length(\tilde f_t)=\tilde\ell_t$, $\cm(\tilde f_t)=\cm( f_t)$,  and $\tilde f_t$ is nonflat. Furthermore, if $\tilde\ell_t/\ell_t\to 1$ as
$t\to 0$, or $1$, then
$
|\tilde f_t-f_t|_1=o(\ell_t).
$
\end{lem}
\begin{proof}
We construct the perturbation $\tilde f_t$  by replacing two
small subarcs, symmetric with respect to the midpoint of $f_t$, by rapidly
oscillating normal bumps, see Figure~\ref{fig:wave}. 
\begin{figure}[h]
\begin{overpic}[height=0.55in]{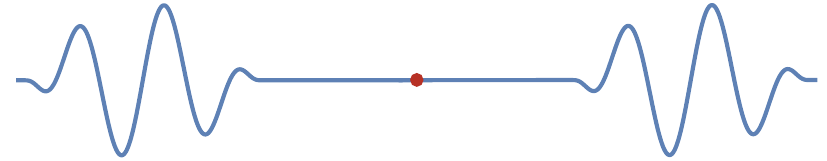}
\end{overpic}
\caption{}\label{fig:wave}
\end{figure}
Their common amplitude controls the increase in
length, while their common distance from the midpoint is adjusted to preserve
the component of the center of mass in the plane of the circle. The reflection
symmetry cancels the components orthogonal to the symmetry plane
automatically; a small additional symmetric perturbation corrects the
remaining component in that plane, and, when $n>3$, further small bumps
ensure that the resulting curve is nonflat.

More precisely, we introduce below a family of general perturbations
$F_{\lambda,q,m,\tau}$, where $\lambda>0$ is a free auxiliary parameter
determining the wavelength, while $q$, $m$, and $\tau$ are to be determined:
$q$ controls the amplitude of the oscillatory bumps, $m$ specifies their
position, and $\tau$ is the size of the additional correction. The argument proceeds in
six stages: (I) we construct the perturbations; (II) we express the
prescribed length and center-of-mass conditions as three equations in
$q,m,\tau$; (III) we let $\lambda\to0$ and obtain a limiting system by
averaging over one period of the oscillation; (IV) we solve this system and
verify that its Jacobian is nonsingular; (V) we apply the implicit function
theorem to solve the exact equations for sufficiently small $\lambda$ and
choose $\lambda$ continuously in $t$; (VI) we establish the estimate at the
endpoints.

\smallskip
\textbf{(I)}
To construct the general perturbation $F_{\lambda,q,m,\tau}$ we
may assume that $I=[-1,1]$. Since the circles have a common center, after
applying a continuous family of rotations fixing that center we may write
$$
f_t(s)=\Bigl(\rho\cos(\sigma_t s),\rho\sin(\sigma_t s),
h,0,\ldots,0\Bigr),
$$
where $\rho=(1-h^2)^{1/2}$, $\rho\sigma_t=\ell_t/2$, and
$0<\sigma_t\leq\pi/4$. For now fix $t$ and suppress it from the notation.
Set $v:=\rho\sigma=\ell/2$ and $r:=\rho\sin(\sigma)/\sigma$, so that
$\cm(f)=(r,0,h,0,\ldots,0)$.
Let $e_1,\ldots,e_n$ be the standard basis of $\R^n$, $R$ be reflection
in the plane spanned by $e_1,e_3$, and set
$$
\nu(s):=
\Bigl(-h\cos(\sigma s),-h\sin(\sigma s),\rho,0,\ldots,0\Bigr).
$$
Then $Rf(s)=f(-s)$, $R\nu(s)=\nu(-s)$, and $\nu$ is a unit normal to
$f$ in $\S^{n-1}$. 

To adjust the first component of the center of mass, we will move the
oscillatory bumps between regions where the first coordinate of $f$ lies above
and below its average value $r$. Choose $d<1/16$ so that, for every
$0<\sigma\leq\pi/4$,
$$
\rho\cos(\sigma s)>r \quad\text{if } |s-\tfrac14|\leq d,
\qquad
\rho\cos(\sigma s)<r \quad\text{if } |s-\tfrac34|\leq d.
$$
Such a uniform choice exists because
$(\cos(\sigma s)-\sin(\sigma)/\sigma)/\sigma^2$ extends continuously to
$\sigma=0$ with value $1/6-s^2/2$, which is positive at $s=1/4$ and
negative at $s=3/4$, while for $\sigma>0$ we have
$\cos(\sigma/4)>\sin(\sigma)/\sigma>\cos(3\sigma/4)$: the first
inequality follows from Vi\`ete's formula
$\sin(\sigma)/\sigma=\prod_{k\geq1}\cos(\sigma/2^k)$, and the second
from $\sin(\sigma)/\sigma\geq1-\sigma^2/6$ and
$\cos(3\sigma/4)\leq1-9\sigma^2/32+27\sigma^4/2048$. Compactness of
$[0,\pi/4]$ then yields the desired $d$.

Fix a nonzero nonnegative function $\alpha\in\C_c^\infty((-d,d))$, which will
determine the oscillatory bumps, and a nonzero nonnegative function
$\psi\in\C_c^\infty((1/32,1/16))$, which will correct the third component
of the center of mass. If $n>3$, choose nonzero functions
$\chi_j\in\C_c^\infty((7/8,15/16))$, $j=4,\ldots,n$, with pairwise
disjoint supports; these will ensure that the resulting curve is nonflat. For $\lambda>0$, $q>0$, $m\in[1/4,3/4]$, and $\tau\in\R$, put, on
$[0,1]$,
$$
u_{\lambda,q,m,\tau}(s):=
\left(
\lambda q \alpha(s-m)\sin\frac{s}{\lambda}+\tau\psi(s)
\right)\nu(s)
+\lambda\sum_{j=4}^n\chi_j(s)e_j,
$$
where the sum is omitted when $n=3$, and define
$$
F_{\lambda,q,m,\tau}(s):=
\frac{f(s)+u_{\lambda,q,m,\tau}(s)}
{|f(s)+u_{\lambda,q,m,\tau}(s)|}.
$$
Extend $F_{\lambda,q,m,\tau}$ to $[-1,1]$ by $F(-s)=RF(s)$.
Note that $u_{\lambda,q,m,\tau}$ vanishes on $[15/16,1]$,  since
$\supp\psi\subset(1/32,1/16)$, the oscillatory bump is supported in
$(m-d,m+d)\subset(3/16,13/16)$ as $m\in[1/4,3/4]$ and $d<1/16$, and
$\supp\chi_j\subset(7/8,15/16)$. Hence, by reflection,
$F_{\lambda,q,m,\tau}=f$ on $[-1,-15/16]\cup[15/16,1]$, the two collars
of width $1/16=|I|/32$. Here
$\lambda$ is the small wavelength which will be chosen later, while
$q,m,\tau$ are the three parameters to be determined.

\smallskip
\textbf{(II)} Now we write  the three equations which $F$ must satisfy.
Since $u_{\lambda,q,m,\tau}$ is normal to $f$, the curve $F$ lies in
$\S^{n-1}$ and agrees with $f$ near $\partial I$. Its speed is even, and
its $e_2,e_4,\ldots,e_n$ components are odd, so the corresponding
components of its center of mass vanish. Set
$$
L(F):=\int_I|F'|\,ds,\quad
C_1(F):=\int_I(F_1-r)|F'|\,ds,\quad
C_3(F):=\int_I(F_3-h)|F'|\,ds.
$$
Thus the  equations to be solved are
\begin{equation}\label{eq:system}
L(F)-\tilde\ell=0,\qquad C_1(F)=0,\qquad C_3(F)=0.
\end{equation}
They are precisely the conditions $\length(F)=\tilde\ell$ and
$\cm(F)=\cm(f)$.

\smallskip
\textbf{(III)}
We next determine the limiting system of equations as $\lambda\to0$, while keeping
$q,m,\tau$ fixed. On the support of the oscillatory bump $\alpha(\,\cdot-m)$, a direct calculation
shows that the speed
$$
|F'_{\lambda,q,m,0}|
=
\sqrt{v^2+q^2\alpha(s-m)^2\cos^2(s/\lambda)}+O(\lambda)
$$
uniformly together with the $q$- and $m$-derivatives. The average of the
leading term over one period of the oscillation is
$\Phi(v,q\alpha(s-m))$, where
$$
\Phi(v,z):=
\frac1{2\pi}\int_0^{2\pi}
\sqrt{v^2+z^2\cos^2\theta}\,d\theta.
$$
Averaging on the support of the oscillatory bump
$\alpha(\,\cdot-m)$, and using ordinary $\C^1$ convergence on the
supports of $\psi$ and the $\chi_j$, therefore show that
$$
(q,m,\tau)\longmapsto
\Bigl(
L(F_{\lambda,q,m,\tau})-\tilde\ell,
C_1(F_{\lambda,q,m,\tau}),
C_3(F_{\lambda,q,m,\tau})
\Bigr)
$$
converges in $\C^1$, locally uniformly in the parameters. Denote the
components of the limiting map by
$\mathcal L(q,m,\tau)$, $\mathcal C_1(q,m,\tau)$, and
$\mathcal C_3(q,m,\tau)$. Thus the limiting equations are
\begin{equation}\label{eq:limit-system}
\mathcal L(q,m,\tau)=0,\qquad
\mathcal C_1(q,m,\tau)=0,\qquad
\mathcal C_3(q,m,\tau)=0.
\end{equation}
When $\tau=0$, their left-hand sides are
\begin{gather*}
\mathcal L(q,m,0)
=
2\int_{-d}^d
\bigl(\Phi(v,q\alpha(u))-v\bigr)\,du
-(\tilde\ell-\ell),\\
\mathcal C_1(q,m,0)
=
2\int_{-d}^d
\bigl(\rho\cos(\sigma(m+u))-r\bigr)
\bigl(\Phi(v,q\alpha(u))-v\bigr)\,du,\\
\mathcal C_3(q,m,0)=0.
\end{gather*}
The asserted convergence follows, for instance, by approximating the
oscillatory integrands and their parameter derivatives by finite Fourier
sums.

\smallskip
\textbf{(IV)}
We now solve the limiting system \eqref{eq:limit-system}. First, with
$\tau=0$, the equation $\mathcal L(q,m,0)=0$ determines the amplitude
$q$. Put
$$
w_q(u):=\Phi(v,q\alpha(u))-v,
\qquad
D(q):=\int_{-d}^d w_q(u)\,du.
$$
The quantity $2D(q)$ is the limiting increase in length produced by the
two symmetric bumps. Since $D$ is strictly increasing from $0$ to infinity,
there is a unique $q_0>0$ such that
$
2D(q_0)=\tilde\ell-\ell.
$
Next, the equation $\mathcal C_1(q_0,m,0)=0$ determines the position
$m$. Its left-hand side is
$$
G(m):=
2\int_{-d}^d
\bigl(\rho\cos(\sigma(m+u))-r\bigr)w_{q_0}(u)\,du.
$$
The choice of $d$ gives $G(1/4)>0>G(3/4)$, while
$$
G'(m)
=
-2\rho\sigma\int_{-d}^d
\sin(\sigma(m+u))w_{q_0}(u)\,du<0.
$$
Hence there is a unique $m_0\in(1/4,3/4)$ such that $G(m_0)=0$.
Finally, $\mathcal C_3(q,m,0)=0$ for all $q,m$, and consequently
$(q_0,m_0,0)$ solves the limiting system. Although no correction is needed
in the limiting system, the parameter $\tau$ will correct the small error
in $C_3$ which occurs when $\lambda>0$.
At $(q_0,m_0,0)$, the derivatives of
$(\mathcal L,\mathcal C_1,\mathcal C_3)$ satisfy
$$
\partial_q\mathcal L=2D'(q_0)>0,\qquad
\partial_m\mathcal C_1=G'(m_0)<0,\qquad
\partial_\tau\mathcal C_3
=
2\rho v\int_0^1\psi(s)\,ds>0.
$$
We also have
$
\partial_m\mathcal L
=
\partial_q\mathcal C_3
=
\partial_m\mathcal C_3
=
0.
$
Thus, with the variables ordered as $(q,m,\tau)$, the determinant of the
Jacobian is
$
\partial_q\mathcal L\,
\partial_m\mathcal C_1\,
\partial_\tau\mathcal C_3\ne0,
$
and hence the Jacobian of the limiting system is invertible.

\smallskip
\textbf{(V)}
We next return to the exact equations \eqref{eq:system}. The implicit function theorem yields,
for every $t$ and all sufficiently small $\lambda>0$, unique parameters
$(q_{t,\lambda},m_{t,\lambda},\tau_{t,\lambda})$ near
$(q_{0,t},m_{0,t},0)$ which solve them. The limiting parameters depend
continuously on $t$, by uniqueness, and the local implicit-function branches
depend continuously on $t$ and agree on overlaps.

For each $t$, let $A(t)>0$ be such that the implicit-function branch is
defined whenever $0<\lambda<A(t)$. By the local uniformity of the preceding
argument, $A$ may be chosen lower semicontinuous. Since every positive lower
semicontinuous function on an interval admits a positive continuous minorant,
we may choose a continuous function $\lambda_t$ with
$0<\lambda_t<A(t)$, and as small as needed below. Then
$
q_t:=q_{t,\lambda_t},
m_t:=m_{t,\lambda_t},
\tau_t:=\tau_{t,\lambda_t}
$
depend continuously on $t$.
Set $\tilde f_t:=F_{\lambda_t,q_t,m_t,\tau_t}$.
Since $u_t$ is orthogonal to $f_t$ and $|f_t|=1$, we have
$$
|\tilde f_t-f_t|_0
\leq
2|u_t|_0
\leq
2\left(
\lambda_tq_t|\alpha|_0
+
|\tau_t||\psi|_0
+
\lambda_t\sum|\chi_j|_0
\right).
$$
The right-hand side tends to zero as $\lambda_t\to0$, locally
uniformly in $t$. Furthermore, for $\lambda$ sufficiently small,
the first two coordinates of
$F_{\lambda,q_{t,\lambda},m_{t,\lambda},\tau_{t,\lambda}}(s)$
are a positive multiple of
$(\cos(\sigma_t s),\sin(\sigma_t s))$. Since
$0<\sigma_t\leq\pi/4$, their polar angle determines $s$, and
hence the curve is embedded. Applying the same
continuous-minorant argument, we may therefore choose
$\lambda_t$ so that $\tilde f_t$ is embedded and
$|\tilde f_t-f_t|_0\leq\epsilon$. It agrees with $f_t$ on the
collars of width $|I|/32$ at $\partial I$, by the vanishing of $u_t$
there recorded in part (I), and has the required length and center
of mass. It is also nonflat: the $\nu$-bump
takes the curve out of the plane of the original circle, and, when $n>3$,
the $\chi_j$-bumps ensure that it is not contained in any hyperplane.

\smallskip
\textbf{(VI)}
It remains to establish the endpoint estimate. Suppose that
$\tilde\ell_t/\ell_t\to1$ at an endpoint. Since $v_t=\ell_t/2$ and
$\Phi(v,z)=v\Phi(1,z/v)$, the equation defining $q_{0,t}$ becomes
$$
\int_{-d}^d
\left[
\Phi\left(1,\frac{q_{0,t}}{v_t}\alpha(u)\right)-1
\right]du
=
\frac{\tilde\ell_t-\ell_t}{\ell_t}.
$$
The left-hand side is a strictly increasing function of $q_{0,t}/v_t$
which vanishes only at zero. Hence $q_{0,t}=o(\ell_t)$.
Applying the same continuous-minorant argument with the additional
requirement
$
|q_{t,\lambda}-q_{0,t}|+|\tau_{t,\lambda}|+\lambda
<\eta(t)\ell_t,
$
where $\eta(t)>0$ is continuous and tends to zero at the endpoint,
we may choose $\lambda_t$ so that
$
|q_t-q_{0,t}|+|\tau_t|+\lambda_t=o(\ell_t).
$
For $u_t:=u_{\lambda_t,q_t,m_t,\tau_t}$ this gives
$|u_t|_0+|u_t'|_0=o(\ell_t)$. Since
$\tilde f_t=(f_t+u_t)/|f_t+u_t|$ and
$|f_t'|_0=\ell_t/2\leq\pi/4$, it follows that
$|\tilde f_t-f_t|_1=o(\ell_t)$.
\end{proof}


\section{Parametric Carath\'{e}odory and Steinitz Theorems}\label{sec:lemmas}
Here we prove one-parameter versions of a pair of classical results in convex geometry.  Steinitz~\cite{steinitz} \cite[Thm.\ 1.3.10]{schneider2014} showed that any point in the interior of the  convex hull of a set $X\subset\R^n$ lies in the interior of the convex hull of at most $2n$ points of $X$. We need the following generalization. We say that a family of sets $S_t=\{s_1(t),\dots, s_N(t)\}\subset I$ is continuous if $s_i(t)$ is continuous.

\begin{lem}\label{lem:steinitz}
Let  $f_t \in \Emb^0(I,\R^{n\geq 3})$ be a $\C^0$-homotopy, and $x \colon [0,1] \to \mathbf{R}^{n}$ be a continuous map with $x(t) \in \inte \conv(f_t)$. Then  there exists a continuous family $S_t\subset \inte(I)$ of sets of $4n$ distinct points  such that $x(t)\in \inte \conv(f_{t}(S_t))$.
\end{lem}
\begin{proof}
By Steinitz's theorem, for each $t_0\in [0,1]$ there is a collection $S_{t_0}\subset I$ of $2n$ distinct points such that $x(t_0)\in \inte\conv(f_{t_0}(S_{t_0}))$. By continuity of $t\mapsto f_t$, for each $t_0\in[0,1]$ we have $x(t)\in \inte \conv(f_{t}(S_{t_0}))$ for $t$ close to $t_0$. By the Lebesgue number lemma, there exists a partition $0=t_0<\dots<t_\ell=1$  such that 
$x(t)\in \inte \conv(f_{t}(S_{t_i}))$ for $t\in[t_i, t_{i+1}]$. Let $S^e_t$ be  a continuous selection of $2n$ points such that 
$S^e_t=S_{t_i}$ whenever $t\in[t_i, t_{i+1}]$  and $i$ is even. Similarly  let $S^o_t$ be  a continuous selection of $2n$ points such that $S^o_t=S_{t_i}$ whenever $t\in[t_i, t_{i+1}]$  and $i$ is odd.
Then $S_t:=S^e_t\cup S^o_t$ gives a continuous selection of $4n$ points which contain $x(t)$ within the interior of their convex hull; however, $S_t$ may contain fewer than $4n$ points. 

To complete the proof it remains to replace $S_t$ by a continuous selection of $4n$ distinct points of $\inte(I)$ which is arbitrarily close to it. Write the coordinates of $S_t$ in nondecreasing order as
$
a\leq\sigma_1(t)\leq\dots\leq\sigma_{4n}(t)\leq b.
$
The order statistics $\sigma_j$ are continuous. For $\eta\in(0,1)$ set
$$
\tilde\sigma_j(t):=(1-\eta)\sigma_j(t)
+\eta\left(a+\frac{j}{4n+1}(b-a)\right).
$$
Then $\tilde\sigma_1(t)<\dots<\tilde\sigma_{4n}(t)$, all these points lie in $\inte(I)$, and the sets $\tilde S_t:=\{\tilde\sigma_1(t),\dots,\tilde\sigma_{4n}(t)\}$ converge uniformly to $S_t$ in Hausdorff distance as $\eta\to0$. Since $x(t)\in\inte\conv(f_t(S_t))$ for every $t$, this inclusion persists under a sufficiently small uniform perturbation, by continuity and compactness of $[0,1]$. Thus, for $\eta$ sufficiently small, $x(t)\in\inte\conv(f_t(\tilde S_t))$ for all $t$.
\end{proof}

Carath\'{e}odory's theorem \cite{caratheodory} \cite[Thm.\ 1.1.4]{schneider2014} states that any point $x$ in the convex hull of a set $X\subset\R^n$ is a convex combination of at most $n+1$ points $p_i$ of $X$, i.e., there exist constants $a_i> 0$ with $\sum a_i=1$ such that $x=\sum a_ip_i$. We need the following continuous version of Carath\'{e}odory's theorem:

\begin{lem}\label{lem:caratheodory}
Let $p_{i}  \colon [0,1]\to\mathbf{R}^{n}$, $i=1,\dots, N$, be $\C^{k}$ maps. Suppose  there exists a  $\C^{k}$ map $x\colon [0,1]\to\mathbf{R}^{n}$ such that $x(t)\in\inte \conv(\{p_i(t)\})$. Then there are $\C^k$ functions $a_{i}(t)$ such that $a_i> 0$, $\sum a_{i}(t) = 1$, and $\sum a_{i}(t) p_{i}(t) =x(t)$.
\end{lem}
\begin{proof}
We may assume that $x(t) =0$ after replacing $p_i(t)$ with $p_{i}(t)-x(t)$.  Then we need to find  functions $a_{i}(t)$, such that 
$
P(t)a(t) = 0,
$
where $P(t)$ is the $n\times N$ matrix with columns $p_i(t)$, and $a(t)$ is the vector with components $a_i(t)$.
So $a(t) \in \ker(P(t))$. Since $0\in\inte\conv\{p_i(t_0)\}$, we may choose $\epsilon\in(0,1/N)$
so small that
$$
-\frac{\epsilon}{1-N\epsilon}\sum p_i(t_0)
\in\conv\big\{p_1(t_0),\dots,p_N(t_0)\big\}.
$$
By Carath\'{e}odory's theorem, we may write this point as   $\sum c_i p_i(t_0)$, where $c_i\geq0$
and $\sum c_i=1$. Setting
$$
a_i^0(t_0):=\epsilon+(1-N\epsilon)c_i,
$$
we obtain $a_i^0(t_0)>0$, $\sum a_i^0(t_0)=1$, and
$\sum a_i^0(t_0)p_i(t_0)=0$.
Note that $P(t)$ has constant rank $n$, since $\inte \conv(\{p_i(t)\})\neq\emptyset$. So $\ker(P(t))$ has dimension $N-n$. Since $p_i(t)$ are $\C^k$, we may choose a $\C^k$ family of orthonormal bases $e_j(t)$ for $\ker(P(t))$ near $t_0$, by a theorem of Dole\v{z}al \cite{dolezal1964,silverman-bucy1970,weiss-falb1969}. Then the projection 
$$
a^0(t):=\sum\left\langle a^0(t_0),e_j(t)\right\rangle e_j(t)
$$ 
of $a^0(t_0)$ into $\ker(P(t))$ will be $\C^k$. Furthermore, since $a_i^0(t_0)>0$,  we have $a_i^0(t)> 0$ for $t$ near $t_0$; and replacing $a^0(t)$ with $a^0(t)/\sum a^0_i(t)$, which again lies in $\ker(P(t))$ since $\ker(P(t))$ is a linear subspace, we may assume that $\sum a^0_i(t)=1$. So the desired coefficients can be found locally. Thus we may cover $[0,1]$ by a finite number of open subsets where the desired coefficients exist and glue them together via a partition of unity.
\end{proof}

\begin{note}
If $p_i(t)$ in the last result were constant, then the conclusion would follow immediately from a theorem of Kalman \cite{kalman1961} who established existence of continuous barycentric coordinates inside convex polytopes.
\end{note}

\begin{note}
Lemma \ref{lem:steinitz} does not hold, for any finite number of points, if we  require that $x(t)$ lie only in the relative interior of $\conv(f_t)$, i.e., allow $f_t$ to be flat \cite{ghomi-mo-2024}.
\end{note}

\section{The Length and Center of Mass of Spherical Curves}\label{sec:spherical}
Here we combine the findings of the last two sections to establish the main
technical result of this work.
\begin{thm}\label{thm:main2}
Let $V\subset\R^n$ be an affine subspace with $\dim(V)\geq 3$ such that $\Sigma:=\S^{n-1}\cap V\neq\emptyset$, $f_t \in\Emb^{k\geq 1}(I,\Sigma)$ be a $\C^1$-isotopy of curves of length $\ell_t$ which are nonflat relative to $V$, and $x\colon[0,1]\to V$ be a continuous map such that $x(t)$ lies in the interior of $\conv(f_t)$ relative to $V$. Suppose that $x(0)=\cm(f_0)$ and $x(1)=\cm(f_1)$.
Then for any sufficiently small neighborhood $W$ of $\partial I$, and any $\epsilon>0$,
 there exists a $\C^1$-isotopy $\tilde f_t\in\Emb^k(I,\S^{n-1})$  such that
\begin{enumerate}
\item{$\tilde f_0=f_0$ and $\tilde f_1=f_1$,}
\item{$|\tilde f_t-f_t|_0\leq\epsilon$,}
\item{$\cm(\tilde f_t)=x(t)$,}
\item{$\tilde f_t=f_t$ on $W$,}
\item{$\tilde f_t$ is nonflat in $\R^{n}$ for $t\in(0,1)$.}
\end{enumerate}
Furthermore, there exists a constant $\mu>0$, depending only
on  $f_t$ and $x(t)$, such that
$\length(\widetilde f_t)$ may be prescribed to equal any
continuous function $\widetilde\ell_t$ satisfying
$\widetilde\ell_0=\ell_0$, $\widetilde\ell_1=\ell_1$,
$\widetilde\ell_t>\ell_t$ for $t\in(0,1)$,
$$
\frac{|x(t)-\cm(f_t)|}{\widetilde\ell_t-\ell_t}
\longrightarrow0
\quad\text{as}\quad t\to0,\,1,
\quad
\text{and}
\quad
\frac{\ell_t|x(t)-\cm(f_t)|}
{\widetilde\ell_t-\ell_t}\leq\mu
\quad\text{on}\quad(0,1).
$$
\end{thm}
\begin{proof}
We first treat the principal case $V=\R^n$, which 
proceeds in four stages: (I) we define a certain point $\ol x(t)$ associated to $x(t)$, and select points $p_i(t)$ of $f_t$ which contain $\ol x(t)$ within their convex hull; (II) we perturb $f_t$ to a curve $\ol f_t$ which contains geodesic segments $\ol g_t^{\,i}$ near $p_i(t)$, see Figure \ref{fig:curves}; (III) we perturb $\ol f_t$ to the desired homotopy $\tilde f_t$ by replacing  $\ol g_t^{\,i}$ with longer curves $\tilde g_t^{\,i}$ via Lemma \ref{lem:bumps}; (IV) we verify that $\tilde f_t$ is  a $\C^1$-homotopy. Finally we address the general case in part (V).

 \begin{figure}[h]
\begin{overpic}[height=0.75in]{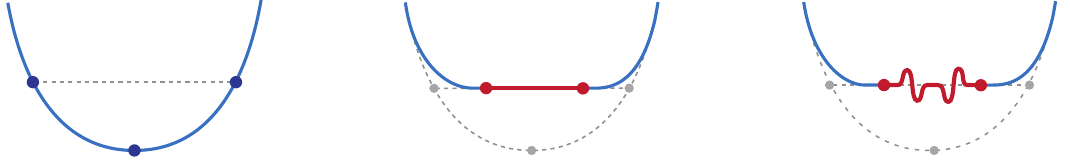}
\put(-2,13){\small$f_t$}
\put(35,13){\small$\ol{f}_t$}
\put(72.5,13){\small$\tilde f_t$}
\put(11,-1){\small$p_i(t)$}
\put(-3.5,7){\small$q_i(t)$}
\put(23.5,7){\small$r_i(t)$}
\put(49,8.5){\small$\ol  g_t^{\,i}$}
\put(86.25,9){\small$\tilde  g_t^{\,i}$}
\end{overpic}
\caption{}\label{fig:curves}
\end{figure}

\textbf{(I)}
Since by assumption $f_t$ is nonflat,
$$
\mu:=\frac12\min_{t\in[0,1]}
\dist\big(x(t),\partial\conv(f_t)\big)>0.
$$
Fix any continuous function $\tilde\ell_t$ satisfying the
conditions of the theorem.
Define $\ol x\colon[0,1]\to\R^n$ by setting
$\ol x(0):=x(0)$, $\ol x(1):=x(1)$, and
$$
\ol x(t):=
\frac{\tilde\ell_t\,x(t)-\ell_t\cm(f_t)}
{\tilde\ell_t-\ell_t}
$$
for $t\in(0,1)$. Since, as  $t\to0,\,1$,
$$
|\ol x(t)-x(t)|
=
\frac{\ell_t\,|x(t)-\cm(f_t)|}
{\tilde\ell_t-\ell_t}
\longrightarrow0
$$
$\ol x$ is continuous. Furthermore,
$$
|\ol x(t)-x(t)|
\leq\mu
<
\dist\big(x(t),\partial\conv(f_t)\big)
$$
for $t\in(0,1)$. Hence
$\ol x(t)\in\inte\conv(f_t)$ for $t\in(0,1)$, while
$\ol x(0)=\cm(f_0)$ and $\ol x(1)=\cm(f_1)$ lie in
$\inte\conv(f_0)$ and $\inte\conv(f_1)$, respectively.
Note that $\dist(\ol x(t),\partial\conv(f_t))\geq\mu$ for all $t$, since
$|\ol x(t)-x(t)|\leq\mu$ and $\dist(x(t),\partial\conv(f_t))\geq2\mu$.
As a collar neighborhood $W$ of $\partial I$ shrinks,
$f_t(I\setminus W)\to f_t(I)$ in Hausdorff distance, uniformly
in $t$ by equicontinuity, and hence so do their convex hulls.
Since
$$
\dist\bigl(\ol x(t),\partial\conv(f_t)\bigr)\geq\mu
$$
for all $t$, it follows that, for every sufficiently small
$W$,
$$
\ol x(t)\in\inte\conv\bigl(f_t(I\setminus W)\bigr)
$$
for every $t$. Fix such a $W$ and set
$I_W:=I\setminus W$.
Now, by Lemma \ref{lem:steinitz}, we may continuously choose
$4n$ distinct points $s_i(t)\in\inte(I_W)$, applying that lemma with
$I_W$ in place of $I$, such that, if
$p_i(t):=f_t(s_i(t))$, then
$$
\ol x(t)\in\inte\conv\{p_i(t)\}.
$$
Thus, by Lemma \ref{lem:caratheodory}, there are continuous
functions $a_i(t)>0$ with $\sum a_i(t)=1$ such that
$$
\sum a_i(t)p_i(t)=\ol x(t).
$$
We assume that all sums range over $i=1,\dots,4n$.

\smallskip
\textbf{(II)}
Let $\delta\colon [0,1]\to\R$ be a continuous function with $\delta(t)>0$ for $t\in(0,1)$ and $\delta(0)=0=\delta(1)$. We may assume that $\max(\delta)$ is smaller than half the smallest distance in $f_t$ between the points $p_i(t)$, $f_t(a)$, $f_t(b)$.
Since $s_i(t)\in\inte(I_W)$ is a continuous family, compactness of
$[0,1]$ allows us to further assume that $\max(\delta)$ is small enough
that the arcs \arc{$r_i(t)q_i(t)$} lie in $f_t(\inte(I_W))$. Then there exist 
points $r_i(t)$, $q_i(t)$ in $f_t$ which are at distance $\delta(t)$ from $p_i(t)$, as measured in $f_t$, and lie on different sides of $p_i(t)$ when $\delta(t)>0$. Let \arc{$r_i(t)q_i(t)$} be the segment of $f_t$ between $r_i(t)$ and $q_i(t)$. So 
$$
|\text{\arc{$r_i(t)q_i(t)$}}|=2\delta(t).
$$
 Replacing \arc{$r_i(t)q_i(t)$}  with the geodesic segment $r_i(t)q_i(t)$ in $\S^{n-1}$, parametrized with constant speed on the parameter interval of the arc it replaces, while keeping the parametrization of $f_t$ elsewhere, yields a $\C^0$-homotopy $\widehat f_t$ of piecewise smooth
embedded curves in $\S^{n-1}$. Since $f_t$ is $\C^1$,  $|\text{\arc{$r_i(t)q_i(t)$}}|/|r_i(t)q_i(t)|\to 1$ as $\delta(t)\to 0$. Thus, if $\max(\delta)$ is sufficiently small, we can make sure that $|r_i(t)q_i(t)|>\delta(t)$. 
 Let $\ol g_t^{\,i}$ be the subsegment of $r_i(t)q_i(t)$ of length $|\ol g_t^{\,i}|=a_i(t)\delta(t)$  centered at the midpoint  
  of $r_i(t)q_i(t)$. So $\ol g_0^{\,i}=p_i(0)$.  If we set $\ol g_t:=\cup _i \ol g_t^{\,i}$, then we may record that
  $$
  |\ol g_t|=\delta(t),\quad\quad\text{and}\quad\quad |\ol g_t^{\,i}|=a_i(t)|\ol g_t|.
  $$
Furthermore, for $t\in(0,1)$, set
$$
\ol p_i(t):=\cm(\ol g_t^{\,i}),
$$
and set
$
\ol p_i(0):=p_i(0),
$
$
\ol p_i(1):=p_i(1).
$
 Next we smooth $\widehat f_t$ by rounding off its corners at
$r_i(t)$ and $q_i(t)$. Extend $\widehat f_t$ slightly beyond $I$. Let $I'\subset\inte(I)$ be a subsegment
so large that $\widehat f_t(I')$ contains all segments
$r_i(t)q_i(t)$, with $I'\subset\inte(I_W)$, $\phi\colon I\to\R$ be a
$\C^\infty$ function such that $0\leq\phi\leq1$, $\phi=1$ on $I'$, and
$\phi=0$ on $W$, and, for $t\in(0,1)$, $\theta_t$ be a continuous
family of symmetric nonnegative $\C^\infty$ mollifiers on $\R$
with $\int_{\R}\theta_t=1$ and sufficiently small compact
support. Define
$$
\ol f_t:=
\frac{\phi(\widehat f_t*\theta_t)+(1-\phi)\widehat f_t}
{\big|\phi(\widehat f_t*\theta_t)+(1-\phi)\widehat f_t\big|}
$$
for $t\in(0,1)$, and set $\ol f_0:=f_0$ and
$\ol f_1:=f_1$.
Since $\ol g_t^{\,i}$ lies in the interior of the
constant-speed geodesic segment $r_i(t)q_i(t)$, we may choose
$\supp(\theta_t)$ so small that convolution there involves only
that segment. By symmetry, $\widehat f_t*\theta_t$ is then a
positive multiple of $\widehat f_t$, so radial projection fixes
$\ol g_t^{\,i}$. Thus, by choosing $\max(\delta)$ and
$\supp(\theta_t)$ sufficiently small, we may assume that
$\ol f_t\in\Emb^k(I,\S^{n-1})$ and
$\ol g_t^{\,i}\subset\ol f_t$.
  
Since the replaced arcs shrink to points and their
tangent jumps tend uniformly to zero as $t\to0$, $1$,
we may choose $\supp(\theta_t)\to0$ so that
$$
|\widehat f_t*\theta_t-f_t|_1\to0.
$$
It follows that $|\ol f_t-f_t|_1\to0$ as
$t\to0$, $1$.
Note also that $\ol f_t=\widehat f_t=f_t$ on $W$, since the replaced
arcs lie in $I_W$, $\phi=0$ on $W$, and $|f_t|=1$; requiring
$\supp(\theta_t)$ small enough that the convolution alters
$\widehat f_t$ only on $\inte(I_W)$ preserves this. The perturbations of part (III) below are supported on parameter
intervals contained in $\inte(I_W)$, so the final curve
$\tilde f_t$ also satisfies $\tilde f_t=f_t$ on $W$, which
yields condition $(4)$.
Furthermore, set
$\ol\ell_t:=\length(\ol f_t)$,
$\widehat\ell_t:=\length(\widehat f_t)$, and
$\epsilon_t:=\delta(t)(\tilde\ell_t-\ell_t)$.
By taking $\supp(\theta_t)$ sufficiently small, we may
arrange that
$$
|\ol\ell_t-\widehat\ell_t|
+
\big|
\ol\ell_t\cm(\ol f_t)
-
\widehat\ell_t\cm(\widehat f_t)
\big|
\leq\epsilon_t
$$
for $t\in(0,1)$, and  consequently
 \begin{equation}\label{eq:l-l}
 |\ol\ell_t-\ell_t|\leq2\big(\ell_t-\widehat\ell_t\big)+\epsilon_t
 =2\sum\big(|\text{\arc{$r_i(t)q_i(t)$}}|-|r_i(t)q_i(t)|\big)+\epsilon_t,
\end{equation}
since $\widehat\ell_t\leq\ell_t$.

\smallskip
\textbf{(III)}
Now we construct the desired homotopy $\tilde f_t$. First note that, for
$t\in(0,1)$,
$$
\ol\ell_t
\leq
\widehat\ell_t+\epsilon_t
\leq
\ell_t+\delta(t)(\tilde\ell_t-\ell_t)
<
\tilde\ell_t,
$$
provided that $\max(\delta)<1$. Define
$$
\hat y(t):=
\frac{\tilde\ell_t x(t)-\ol\ell_t\cm(\ol f_t)}
{\tilde\ell_t-\ol\ell_t}.
$$
This is the analogue of $\ol x(t)$ with $\ol f_t$ in place of $f_t$. A
direct computation gives
$$
\hat y(t)-\ol x(t)
=
\frac{
\ell_t\cm(f_t)-\ol\ell_t\cm(\ol f_t)
+(\ol\ell_t-\ell_t)\ol x(t)
}{
\tilde\ell_t-\ol\ell_t
}.
$$
Thus $\hat y(t)\to\ol x(t)$ as $t\to0$, $1$, provided that
$$
\frac{|\ol\ell_t-\ell_t|}
{\tilde\ell_t-\ell_t}\to0
\qquad\text{and}\qquad
\frac{
\big|
\ol\ell_t\cm(\ol f_t)-\ell_t\cm(f_t)
\big|
}{
\tilde\ell_t-\ell_t
}
\to0.
$$
These limits will be arranged in part (IV) through the choice of
$\delta(t)$. Since
$\ol x(t)\in\inte\conv\{p_i(t)\}$ for every $t$, compactness
implies that this inclusion persists under sufficiently small
uniform perturbations of $\ol x(t)$ and  $p_i(t)$.
Near $t=0$, $1$, the required perturbations tend to zero by the
limits above and the fact that $\ol p_i(t)\to p_i(t)$. On every
compact subinterval of $(0,1)$, they may be made uniformly small
by choosing $\max(\delta)$ and $\supp(\theta_t)$ sufficiently
small, since $\tilde\ell_t-\ell_t$ is bounded away from zero
there. Thus we may assume that
$
\hat y(t)\in\inte\conv\{\ol p_i(t)\}
$
for $t\in(0,1)$. Set $\hat y(0):=x(0)$ and
$\hat y(1):=x(1)$. Then $\hat y$ is continuous on $[0,1]$, and
the same inclusion holds at the endpoints.  Hence, by Lemma
\ref{lem:caratheodory}, there are continuous functions $b_i(t)>0$ such
that
$
\sum b_i(t)=1,
$
$
\sum b_i(t)\ol p_i(t)=\hat y(t).
$
Prescribe the lengths of the replacement curves by
$$
|\tilde g_t^{\,i}|
:=
|\ol g_t^{\,i}|
+
b_i(t)(\tilde\ell_t-\ol\ell_t).
$$
Then $|\tilde g_t^{\,i}|>|\ol g_t^{\,i}|$ for $t\in(0,1)$, and
$
\sum|\tilde g_t^{\,i}|
=
|\ol g_t|+\tilde\ell_t-\ol\ell_t.
$
The segments $\ol g_t^{\,i}$, being geodesics of length at most
$\max(\delta)$, subtend angles at most $\pi/2$ if
$\max(\delta)\leq\pi/2$. Lemma \ref{lem:bumps} therefore supplies curves
$\tilde g_t^{\,i}$ of the prescribed lengths such that
$$
\cm(\tilde g_t^{\,i})
=
\cm(\ol g_t^{\,i})
=
\ol p_i(t).
$$
Let $\tilde f_t$ be obtained from $\ol f_t$ by making these replacements,
and set $\tilde g_t:=\bigcup_i\tilde g_t^{\,i}$. Then
$$
\length(\tilde f_t)
=
\ol\ell_t-|\ol g_t|
+\sum|\tilde g_t^{\,i}|
=
\tilde\ell_t.
$$
Furthermore,
\begin{gather*}
\tilde\ell_t\cm(\tilde f_t)
=
\ol\ell_t\cm(\ol f_t)
+
\sum
\bigl(
|\tilde g_t^{\,i}|-|\ol g_t^{\,i}|
\bigr)\ol p_i(t)
=
\ol\ell_t\cm(\ol f_t)
+
(\tilde\ell_t-\ol\ell_t)
\sum b_i(t)\ol p_i(t)\\
=
\ol\ell_t\cm(\ol f_t)
+
(\tilde\ell_t-\ol\ell_t)\hat y(t)
=
\tilde\ell_t x(t).
\end{gather*}
Thus $\cm(\tilde f_t)=x(t)$. By Lemma \ref{lem:bumps}, the replacement
curves, and consequently $\tilde f_t$, may also be required to be nonflat
in $\R^{n}$ for $t\in(0,1)$.

We also record that condition (2) of Theorem \ref{thm:main2} holds. On the parameter interval of each replaced arc, $f_t(s)$ and $\widehat f_t(s)$ both lie within $2\delta(t)$ of $r_i(t)$, since the arc and the chord each have length at most $2\delta(t)$; hence $|\widehat f_t-f_t|_0\leq4\max(\delta)$. Choosing $\supp(\theta_t)$ so small that $|\ol f_t-\widehat f_t|_0\leq\max(\delta)$, and applying Lemma \ref{lem:bumps} with $\epsilon/2$ in place of $\epsilon$, we obtain $|\tilde f_t-f_t|_0\leq5\max(\delta)+\epsilon/2\leq\epsilon$, for $\max(\delta)\leq\epsilon/10$.
 
\smallskip 
\textbf{(IV)}
It remains to check that $\tilde f_t$ is a $\C^1$-homotopy, i.e.,  $|\tilde f_t-f_0|_1\to 0$  as $t\to 0$, and likewise for $f_1$ as $t\to 1$; we consider the first limit, as the second one is analogous. To this end note that
$
|\tilde f_t-f_0|_1 \leq |\tilde f_t-\ol f_t|_1 + |\ol f_t-f_0|_1
$
and $|\ol f_t-f_0|_1\leq |\ol f_t-f_t|_1+|f_t-f_0|_1\to 0$, as we discussed in Part II above, and since $f_t$ is a $\C^1$-homotopy. So we just need to show that $|\tilde f_t-\ol f_t|_1\to 0$, which is the case provided that $|\tilde  g^{\,i}_t-\ol g^{\,i}_t|_1\to 0$. For each $i$, let $J_t^i\subset I$ be the parameter
interval corresponding to $\ol g_t^{\,i}$, and let
$\psi_t^i\colon[-1,1]\to J_t^i$ be the
orientation-preserving affine map. The application of
Lemma \ref{lem:bumps} in part (III) is understood after
pulling $\ol g_t^{\,i}$ back to $[-1,1]$ by $\psi_t^i$.
Thus, if
$
G_t^i:=\ol f_t\circ\psi_t^i
$
and $\tilde G_t^i$ is the curve supplied by the lemma,
then
$
|\tilde G_t^i-G_t^i|_1=o(|\ol g_t^{\,i}|)
$
provided that
$
|\tilde g_t^{\,i}|/|\ol g_t^{\,i}|\to1.
$
Since $\ol f_t\to f_0$ in $\C^1$ and $f_0$ is an
immersion, $|J_t^i|$ and $|\ol g_t^{\,i}|$ are
comparable for $t$ close to $0$. Consequently,
$$
|\tilde g_t^{\,i}-\ol g_t^{\,i}|_0
=
|\tilde G_t^i-G_t^i|_0
\longrightarrow0,
$$
$$
|(\tilde g_t^{\,i})'-(\ol g_t^{\,i})'|_0
=
\frac{2}{|J_t^i|}
|(\tilde G_t^i)'-(G_t^i)'|_0
\longrightarrow0.
$$
It therefore suffices to check that
$
|\tilde g_t^{\,i}|/|\ol g_t^{\,i}|\to1
$
as $t\to0$, $1$. Note that, for $t\in (0,1)$,
$$
\frac{|\tilde  g^{\,i}_t|}{|\ol g^{\,i}_t|}
=
1+\frac{b_i(t)}{a_i(t)}\cdot\frac{\tilde\ell_t-\ol\ell_t}{\delta(t)},
$$
where $b_i/a_i$ is bounded on $[0,1]$, since $a_i$, $b_i$ are continuous and positive. So it remains to check that $(\tilde\ell_t-\ol\ell_t)/\delta(t)\to 0$ as $t\to 0$, $1$.
To this end, note that 
$$
\frac{\tilde\ell_t-\ol\ell_t}{\delta(t)}\leq \frac{|\tilde\ell_t-\ell_t|}{\delta(t)}+\frac{|\ell_t-\ol\ell_t|}{\delta(t)}.
$$
We can make sure that
$(\tilde\ell_t-\ell_t)/\delta(t)$ vanishes as
$t\to0$, $1$, and at the same time obtain both limits
required in part (III). Indeed, for constant $\delta'$,
let $\widehat f_t(\delta')$ and
$\widehat\ell_t(\delta')$ denote the corresponding curve
and its length, and set
$$
\omega(\delta):=
\delta+
\sup_{t\in[0,1],\;0<\delta'\leq\delta}\;
\frac{2}{\delta'}
\left(
\ell_t-\widehat\ell_t(\delta')
+
\big|
\widehat\ell_t(\delta')\cm(\widehat f_t(\delta'))
-
\ell_t\cm(f_t)
\big|
\right).
$$
Then $\omega(\delta)\to0$ as $\delta\to0$. Replacing
$\omega$, if necessary, by a continuous increasing function
which is everywhere at least $\omega$ and still tends to
zero at zero, choose $\delta(t)$, for $t$ close to
$0$, $1$, so that
$$
\omega\big(\delta(t)\big)^{1/2}\delta(t)
=
\tilde\ell_t-\ell_t.
$$
Since $\omega(\delta)\geq\delta$, the function
$\delta\mapsto\omega(\delta)^{1/2}\delta$ is continuous
and strictly increasing, so this determines
$\delta(t)$ continuously. 
After shrinking the endpoint neighborhoods if necessary, extend
$\delta$ continuously and positively over the remaining compact
subinterval, keeping $\max(\delta)$ sufficiently small to satisfy
all the requirements imposed in parts (II) and (III).
By the
definition of $\omega$ and the choice made in part (II),
both ratios in part (III) are bounded by
$
\omega(\delta(t))^{1/2}+\delta(t),
$
which tends to zero as $t\to0$, $1$. Furthermore,
$$
\frac{\tilde\ell_t-\ell_t}{\delta(t)}
=
\omega\big(\delta(t)\big)^{1/2}\longrightarrow0.
$$
Since the first ratio required in part (III) also tends to zero,
we have
$$
\frac{\tilde\ell_t-\ol\ell_t}{\delta(t)}
\leq
\frac{\tilde\ell_t-\ell_t}{\delta(t)}
\left(
1+
\frac{|\ol\ell_t-\ell_t|}
{\tilde\ell_t-\ell_t}
\right)
\longrightarrow0.
$$
Hence
$|\tilde g_t^{\,i}|/|\ol g_t^{\,i}|\to1$, as desired.

Finally note that $\tilde f_t$ is embedded. Near $t=0$, $1$, the estimates above and the last assertion of Lemma \ref{lem:bumps} show that $\tilde f_t$ is $\C^1$-close to $f_t$, and hence is embedded. On any compact subinterval of $(0,1)$, the portions of the segments $\ol g_t^{\,i}$ on which the perturbations are supported are uniformly separated from one another and from $\ol f_t\setminus\ol g_t$. Fix $\eta_0>0$ so that $\tilde f_t$ is embedded for $t\in(0,\eta_0]\cup[1-\eta_0,1)$ by the preceding estimates, which do not involve the closeness parameter of Lemma \ref{lem:bumps}. Then, in applying that lemma, we choose its closeness parameter smaller than both $\epsilon/2$, as in part (III), and the separation  corresponding to $[\eta_0,1-\eta_0]$, so that these perturbations remain disjoint from one another and from the rest of $\ol f_t$ for $t\in[\eta_0,1-\eta_0]$. Since each $\tilde g_t^{\,i}$ is embedded and agrees with $\ol g_t^{\,i}$ near its endpoints, $\tilde f_t$ is embedded there as well. This completes the proof in the case $V=\R^n$.

\smallskip 
\textbf{(V)}
Finally, suppose that $V$ is a proper affine subspace of $\R^n$. Then $\Sigma$ is a sphere of radius $\rho:=(1-|q_0|^2)^{1/2}\leq 1$ centered at the point $q_0$ of $V$ nearest the origin. The above argument applies with $\Sigma$ in place of $\S^{n-1}$, as follows. Stages (I) and (II) are carried out within $\Sigma$: Lemmas \ref{lem:steinitz} and \ref{lem:caratheodory} apply in $V$, since $\dim(V)\geq 3$ and $f_t$ is nonflat relative to $V$; the constant $\mu$ and the distance to the boundary of $\conv(f_t)$ in stage (I) are taken relative to $V$; the segments $r_i(t)q_i(t)$ are taken to be geodesics of $\Sigma$; and the normalization in the definition of $\ol f_t$ is replaced by the radial projection 
$$
z\longmapsto q_0+\rho(z-q_0)/|z-q_0|
$$ 
onto $\Sigma$, which fixes the geodesics of $\Sigma$ pointwise on the corresponding segments, as before. In particular $\ol f_t\subset\Sigma$, and the segments $\ol g_t^{\,i}$, being geodesics of $\Sigma$, are circular arcs of $\S^{n-1}$, which subtend angles at most $\pi/2$ in their circles, assuming $\max(\delta)\leq\pi\rho/2$. Thus, in stage (III), Lemma \ref{lem:bumps} yields curves $\tilde g_t^{\,i}$ in $\S^{n-1}$ with $\cm(\tilde g_t^{\,i})=\ol p_i(t)$ and prescribed lengths; since the length and center-of-mass computations above use only these prescribed lengths and centers of mass, they are unchanged. Furthermore, $\tilde f_t$ is nonflat in $\R^{n}$ for $t\in(0,1)$, since so are the curves $\tilde g_t^{\,i}$ by Lemma \ref{lem:bumps}. Stage (IV) is unchanged, which completes the proof of Theorem \ref{thm:main2}.
\end{proof}

\section{Proof of the Main Result} \label{sec:proof}
Here we establish Theorem \ref{thm:main1}. We need only one more lemma, which produces arbitrarily small deformations through curves of constant curvature whose tantrices are nonflat on each
subinterval of a prescribed partition. 

\begin{lem}\label{lem:endpoints}
Let $f_0\in\Imm^{k\geq 2}(\Gamma,\R^{n\geq 3})$ be a unit speed curve with constant curvature $\kappa_0>0$, $\epsilon>0$, and $I_i$ be subintervals in a partition of $\Gamma$ such that each restriction $T_0^i$ of the tantrix $T_0$ of $f_0$ to $I_i$ is embedded. Then there exists a $\C^2$-homotopy $\tilde f_t\in\Imm^k(\Gamma,\R^n)$ of unit speed curves with constant curvature such that $\tilde f_0=f_0$, $|\tilde f_t-f_0|_1\leq\epsilon$, and, for $t\in (0,1)$, the tantrix $\tilde T_t$ of $\tilde f_t$ is nonflat on each $I_i$ and $\tilde T_t|_{I_i}$ is embedded.
\end{lem}
\begin{proof}
Since $f_0$ has unit speed, $\kappa_0:=|T_0'|$ is the constant curvature of $f_0$, and $T_0$ has constant speed $\kappa_0$. Set $T_t:=T_0$ and $\epsilon':=\epsilon/(2(|\Gamma|+1))$, where $|\Gamma|:=b-a$ is the length of the parameter domain.
Refine each $I_i$ into finitely many subintervals $I_{ij}$ so that
$
\operatorname{diam}(T_0(I_{ij}))<\epsilon'.
$
If $T_0|_{I_{ij}}$ is a circular arc which subtends an angle greater than $\pi/2$, quadrisect $I_{ij}$. Since an embedded circular arc subtends an angle less than $2\pi$, we may therefore assume that every circular restriction $T_0^{ij}:=T_0|_{I_{ij}}$ subtends an angle at most $\pi/2$.

Let $c\colon[0,1]\to\R$ be a continuous function with $c(0)=c(1)=\kappa_0$ and $c(t)>\kappa_0$ for $t\in(0,1)$. Set $T_t^{ij}:=T_0^{ij}$, $x_{ij}(t):=\ave(T_0^{ij})$, and
$
\tilde\ell_t^{ij}:=c(t)|I_{ij}|.
$
The endpoint and strict inequality requirements for $\tilde\ell_t^{ij}$ follow from the choice of $c$, while the two displayed conditions in the last statement of Theorem \ref{thm:main2} hold trivially, since $x_{ij}(t)=\cm(T_t^{ij})$ for all $t$.

Let $V_{ij}$ be the affine hull of $T_0(I_{ij})$. Then $T_0^{ij}$ is
nonflat relative to $V_{ij}$, and $\dim(V_{ij})\geq2$. If $\dim(V_{ij})\geq3$, choose a sufficiently small collar
neighborhood $W_{ij}$ of $\partial I_{ij}$ as allowed by item
$(4)$ of Theorem \ref{thm:main2}, let $\eta_{ij}$ be its width,
and apply that theorem to $T_t^{ij}$, with respect to
$x_{ij}(t)$ and $\tilde\ell_t^{ij}$, with $V:=V_{ij}$,
$W:=W_{ij}$, and closeness parameter $\epsilon_{ij}>0$ to be
chosen below. If $\dim(V_{ij})=2$, apply Lemma
\ref{lem:bumps} with the same closeness parameter to the
constant family $T_t^{ij}$, for $t\in(0,1)$, with respect to
$\tilde\ell_t^{ij}$, set
$\tilde T_0^{ij}:=T_0^{ij}=:\tilde T_1^{ij}$, and set
$\eta_{ij}:=|I_{ij}|/32$. The last statement of that lemma
gives continuity at the endpoints, since
$c(t)/\kappa_0\to1$.
In either case, the resulting curve $\tilde T_t^{ij}$ is embedded,
nonflat in $\R^n$ for $t\in(0,1)$, agrees with $T_0$ near
$\partial I_{ij}$, and satisfies
$$
\cm(\tilde T_t^{ij})=\ave(T_0^{ij}),
\quad
\length(\tilde T_t^{ij})=c(t)|I_{ij}|,
\quad
|\tilde T_t^{ij}-T_0^{ij}|_0\leq\epsilon_{ij},
\quad
\tilde T_t^{ij}=T_0^{ij},
$$
on collars of width $\eta_{ij}$.
Set $\eta:=\min_j\eta_{ij}$. Since $T_0|_{I_i}$ is a continuous
injection of a compact interval, there is $\beta>0$ such that
$|T_0(s)-T_0(s')|\geq\beta$ whenever $s$, $s'\in I_i$ and
$|s-s'|\geq\eta$. Choose $\epsilon_{ij}<\min\{\epsilon',\beta/4\}$. Then
the pieces over $I_i$ glue to an injective map: if
$\tilde T_t^{ij}(s)=\tilde T_t^{ij'}(s')$ with $j\neq j'$, then
$|T_0(s)-T_0(s')|\leq2\max_j\epsilon_{ij}<\beta$, so $|s-s'|<\eta$,
which places $s$ and $s'$ in the collars of their pieces, where these
coincide with $T_0$; hence $s=s'$ is a shared subdivision point, where
adjacent pieces agree. As this argument involves only the images of the
pieces and their values at subdivision points, it persists under the
reparametrizations below.
Parametrize each $\tilde T_t^{ij}$ with constant speed $c(t)$ and glue
them to obtain $\tilde T_t$ on $\Gamma$. Near every subdivision point, the adjacent pieces are
orientation-preserving constant-speed reparametrizations of the same
portion of $T_0$, with common speed $c(t)$;
hence they glue with class $\C^{k-1}$, and
$\tilde T_t|_{I_i}$ is embedded. Since every $\tilde T_t^{ij}$ is
nonflat, $\tilde T_t$ is nonflat on $I_i$ for $t\in(0,1)$.
Furthermore, since $\tilde T_t^{ij}$ has constant speed,
$$
\int_{I_{ij}}\tilde T_t
=
|I_{ij}|\cm(\tilde T_t^{ij})
=
|I_{ij}|\ave(T_0^{ij})
=
\int_{I_{ij}}T_0.
$$
Identify $\Gamma$ with $I=[a,b]$, taking $a$ to be a subdivision point, and set
$$
\tilde f_t(s):=f_0(a)+\int_a^s\tilde T_t.
$$
Then $\tilde f_t$ has unit speed and constant curvature $c(t)$, agrees with $f_0$ at all subdivision points, and is closed and $\C^k$ at $a$ when $\Gamma$ is a circle. Also, if $s\in I_{ij}$, then $\tilde T_t(s)$ belongs to the
image of $\tilde T_t^{ij}$. Hence there is $s'\in I_{ij}$ such
that
$$
|\tilde T_t(s)-T_0(s')|\leq\epsilon_{ij}.
$$
Since $\operatorname{diam}(T_0(I_{ij}))<\epsilon'$, it follows
that
$
|\tilde T_t(s)-T_0(s)|
<
\epsilon_{ij}+\epsilon'
<
2\epsilon'.
$
Consequently,
$$
|\tilde f_t-f_0|_1
\leq
(|\Gamma|+1)|\tilde T_t-T_0|_0
\leq
2\epsilon'(|\Gamma|+1)
=
\epsilon.
$$
Finally, $\tilde f_0=f_0=\tilde f_1$.
\end{proof}

Finally we are ready to establish the main result of this work: 

\begin{proof}[Proof of Theorem \ref{thm:main1}]
The argument proceeds in five stages: (I) we normalize
the homotopy and choose a common partition; (II) we arrange that the
tantrix is nonflat on each member of the partition; (III) we choose the
target curvature and verify the length hypotheses of Theorem
\ref{thm:main2}; (IV) we apply Theorem
\ref{thm:main2} to the individual tantrix
arcs and glue the resulting curves; (V) we integrate the new tantrices and
undo the preceding reductions.

\smallskip
\textbf{(I)}
We first normalize the homotopy and choose the partition. Identify
$\Gamma$ with $I=[a,b]$, with the endpoints identified when $\Gamma$ is
a circle. Set $c_t:=\length(f_t)/|I|$ and replace $f_t$ by
$c_t^{-1}f_t$. The resulting curves have constant length $|I|$, and
there is a continuous family of diffeomorphisms
$\theta_t\colon I\to I$ such that $c_t^{-1}f_t\circ\theta_t$ has unit
speed. If $\widehat f_t$ is a homotopy obtained for this normalized
family, then
$$
\tilde f_t:=c_t\widehat f_t\circ\theta_t^{-1}
$$
is a homotopy for the original family. Moreover, sufficiently close
approximation before undoing these normalizations gives the required
approximation afterward; see the proof of
\cite[Prop.\ 5.4]{ghomi:knots}. We may therefore assume that $f_t$ has
unit speed and constant length $|I|$. These normalizations preserve the
constancy of the curvatures of $f_0$ and $f_1$.
Let $T_t:=f_t'$. Fix $\rho>0$ so small that
$
2(|I|+1)\rho<\epsilon/4.
$
After subdividing $I$ into finitely many subintervals $I_i$, we may
assume that each $T_t|_{I_i}$ is embedded and
$
\operatorname{diam}(T_t(I_i))<\rho/2
$
for every $t$ and $i$.

\smallskip
\textbf{(II)}
We next arrange that $T_t|_{I_i}$ is nonflat for every $t$ and $i$.
Apply Lemma \ref{lem:endpoints} to $f_0$ and $f_1$, with respect to the
partition chosen in (I) and with
$\min\{\epsilon/8,\rho/4\}$ in place of $\epsilon$. Let $g_t$ and $h_t$
be the resulting homotopies, so $g_0=f_0$ and $h_0=f_1$.
Fix $t^*\in(0,1)$, and set
$\ol f_0:=g_{t^*}$ and $\ol f_1:=h_{t^*}$. Their tantrices are nonflat
on each $I_i$. Choose $\eta\in(0,1/4)$ sufficiently small that
$$
\sup_{t,t'\in[0,1],\; |t-t'|\leq2\eta}
|f_t-f_{t'}|_1
<
\frac{\epsilon}{2}.
$$
For $u\in[0,1]$, let $F_u^0$ be the concatenation of the reverse of
$g|_{[0,t^*]}$, traversed on $[0,\eta]$, the family
$f_{\phi(u)}$ on $[\eta,1-\eta]$, where
$
\phi(u):=(u-\eta)/(1-2\eta),
$
and $h|_{[0,t^*]}$, traversed on $[1-\eta,1]$. Thus $F_u^0$ is a
homotopy of unit speed curves of constant length from $\ol f_0$ to
$\ol f_1$, and its endpoint tantrices are nonflat on every $I_i$. Apply Proposition \ref{prop:nonflat} to the restrictions
$F_u^0|_{I_i}$. Choosing the perturbations sufficiently small, they glue
to a homotopy $F_u$ such that
$$
|F_u-F_u^0|_1<\frac{\epsilon}{8},
$$
each restriction of its tantrix to $I_i$ is embedded and nonflat, and the
diameter bound in (I) remains valid. The perturbations preserve
the length on each $I_i$ and agree with $F_u^0$ near its endpoints;
hence, after reparametrizing each restriction with unit speed, the pieces
still glue smoothly. Relabeling $u$ as $t$, let $T_t$ and $\kappa_t$
denote the tantrix and curvature of $F_t$. Thus $F_t$ has unit speed and
constant length, $T_t|_{I_i}$ is embedded and nonflat for every $t$ and
$i$, $\operatorname{diam}(T_t(I_i))<\rho$, and $F_0$, $F_1$ have
constant curvature.

\smallskip
\textbf{(III)}
We now choose the curvature of the desired homotopy and verify the
hypotheses needed to prescribe the lengths of the tantrix arcs. Set
$$
\omega(t):=\max_I\kappa_t-\min_I\kappa_t,
\qquad
\kappa_-:=\min_{t,s}\kappa_t(s)>0,
\qquad
\kappa_+:=\max_{t,s}\kappa_t(s).
$$
Since $F_0$ and $F_1$ have constant curvature,
$\omega(t)\to0$ as $t\to0$, $1$. Choose a continuous function
$c\colon[0,1]\to\R$ whose endpoint values are the curvatures of
$F_0$, $F_1$, and such that
$$
c(t)>\max_I\kappa_t+K\sqrt{\omega(t)}
$$
for $t\in(0,1)$, where $K>0$ will be chosen below.
For $T_t^i:=T_t|_{I_i}$, set
$$
x_i(t):=\ave(T_t^i),
\qquad
\ell_t^i:=\length(T_t^i),
\qquad
\tilde\ell_t^i:=c(t)|I_i|.
$$
Since $T_t^i$ is nonflat,
$x_i(t)\in\inte\conv(T_t^i)$. Furthermore,
$x_i(0)=\cm(T_0^i)$ and $x_i(1)=\cm(T_1^i)$, because the endpoint
curvatures are constant. We also have
$\tilde\ell_0^i=\ell_0^i$, $\tilde\ell_1^i=\ell_1^i$, and
$\tilde\ell_t^i>\ell_t^i$ for $t\in(0,1)$. Since
$|T_t'|=\kappa_t$,
$$
\cm(T_t^i)-x_i(t)
=
\int_{I_i}T_t
\left(
\frac{\kappa_t}{\int_{I_i}\kappa_t}
-
\frac{1}{|I_i|}
\right)ds,
$$
and consequently
$$
|x_i(t)-\cm(T_t^i)|
\leq
\frac{\omega(t)}{\kappa_-},
\qquad
\tilde\ell_t^i-\ell_t^i
\geq
K\sqrt{\omega(t)}\,|I_i|.
$$
It follows that
$$
\frac{|x_i(t)-\cm(T_t^i)|}
{\tilde\ell_t^i-\ell_t^i}
\leq
\frac{\sqrt{\omega(t)}}{K\kappa_-|I_i|}
\longrightarrow0
$$
as $t\to0$, $1$. Also, since
$\ell_t^i\leq\kappa_+|I_i|$,
$$
\frac{\ell_t^i|x_i(t)-\cm(T_t^i)|}
{\tilde\ell_t^i-\ell_t^i}
\leq
\frac{\kappa_+\sqrt{\omega(t)}}{K\kappa_-}.
$$
Let $\mu$ be the smallest of the constants in the last statement of
Theorem \ref{thm:main2} for the finitely many families
$T_t^i$, $x_i(t)$. Choosing $K$ sufficiently large, we obtain the
required bound by $\mu$. Thus all the hypotheses needed to prescribe the
lengths $\tilde\ell_t^i$ are satisfied.

\smallskip
\textbf{(IV)}
Apply Theorem \ref{thm:main2} to each family $T_t^i$, with
$V=\R^n$, prescribed center of mass $x_i(t)$, prescribed length
$\tilde\ell_t^i$, $\rho$ as the closeness parameter, and a
sufficiently small neighborhood of $\partial I_i$. We obtain curves
$\tilde T_t^i$ such that
$$
\cm(\tilde T_t^i)=\ave(T_t^i),
\qquad
\length(\tilde T_t^i)=c(t)|I_i|,
\qquad
|\tilde T_t^i-T_t^i|_0\leq\rho.
$$
Parametrize each $\tilde T_t^i$ with constant speed $c(t)$ and define
$\tilde T_t$ by gluing these curves along the partition. The pieces
glue to a $\C^{k-1}$ curve, since near every partition point the
adjacent pieces are constant-speed, orientation-preserving
reparametrizations of the same curve $T_t$. Moreover, if $s\in I_i$,
then $\tilde T_t(s)$ belongs to the image of $\tilde T_t^i$, and hence
there is $s'\in I_i$ such that
$
|\tilde T_t(s)-T_t(s')|\leq\rho.
$
Since $\operatorname{diam}(T_t(I_i))<\rho$, it follows that
$
|\tilde T_t(s)-T_t(s)|<2\rho.
$

\smallskip
\textbf{(V)}
Define
$
\widehat F_t(s):=F_t(a)+\int_a^s\tilde T_t.
$
Then $\widehat F_t$ has unit speed and constant curvature $c(t)$. Since
$\tilde T_t^i$ has constant speed,
$
\ave(\tilde T_t^i)
=
\cm(\tilde T_t^i)
=
\ave(T_t^i),
$
and hence
$$
\int_{I_i}\tilde T_t
=
\int_{I_i}T_t.
$$
Thus $\widehat F_t$ agrees with $F_t$ at all partition points. When
$\Gamma$ is a circle, it is closed and $\C^k$ at the identified endpoint
by the same gluing argument. The endpoint conditions in Theorem
\ref{thm:main2} also give
$\widehat F_0=F_0$ and $\widehat F_1=F_1$. Furthermore,
$$
|\widehat F_t-F_t|_1
\leq
(|I|+1)|\tilde T_t-T_t|_0
<
2(|I|+1)\rho
<
\frac{\epsilon}{4}.
$$
Concatenate $g|_{[0,t^*]}$, the homotopy $\widehat F_t$, and the reverse
of $h|_{[0,t^*]}$, traversing the first and last pieces on $[0,\eta]$
and $[1-\eta,1]$, respectively. The resulting normalized homotopy,
denoted by $\widehat f_t$, joins $f_0$ to $f_1$ through curves of
constant curvature. The curvature varies continuously at the junctions
because $c(0)$ and $c(1)$ are the curvatures of $\ol f_0$ and
$\ol f_1$. By the choice of $\eta$ in (II),
$$
\sup_{t,t'\in[0,1],\; |t-t'|\leq2\eta}
|f_t-f_{t'}|_1
<
\frac{\epsilon}{2}.
$$Together with the choices in (II) and the estimate above, this
yields
$$
|\widehat f_t-f_t|_1
<
\frac{\epsilon}{8}
+
\frac{\epsilon}{8}
+
\frac{\epsilon}{4}
+
\frac{\epsilon}{2}
=
\epsilon.
$$
Finally, undoing the rescaling and reparametrization from (I)
yields the required homotopy for the original family.
\end{proof}

\bibliographystyle{amsplain}
\bibliography{references}

\end{document}